\documentclass[12pt]{article}
\usepackage[english]{babel}
\usepackage{amssymb}
\usepackage{amsmath}
\usepackage{mathrsfs}
\newtheorem{thm}{Theorem}

\begin{document}

\begin{center}
\large{\bf FINITENESS OF HITTING TIMES UNDER TABOO}
\end{center}
\vskip0,5cm
\begin{center}
Ekaterina Vl. Bulinskaya\footnote{ \emph{Email address:} {\tt
bulinskaya@yandex.ru}}$^,$\footnote{The work is partially supported
by Dmitry Zimin Foundation ``Dynasty''.} \vskip0,2cm \emph{Lomonosov
Moscow State University}
\end{center}
\vskip1cm

\begin{abstract}
We consider a continuous-time Markov chain with a finite or
countable state space. For a site $y$ and subset $H$ of the state
space, the hitting time of $y$ under taboo $H$ is defined to be
infinite if the process trajectory hits $H$ before $y$, and the
first hitting time of $y$ otherwise. We investigate the probability
that such times are finite. In particular, if the taboo set is
finite, an efficient iterative scheme reduces the study to the known
case of a singleton taboo. A similar procedure applies in the case
of finite complement of the taboo set. The study is motivated by
classification of branching processes with finitely many catalysts.

\vskip0,5cm {\it Keywords and phrases}: Markov chain, hitting time,
taboo probabilities, catalytic branching process.

\vskip0,5cm 2010 {\it AMS classification}: 60J27, 60G40.
\end{abstract}

\section{Introduction}

The concept of passage time under taboo for Markov chain has a long
history and the first comprehensive exposition of the subject was
given in the classical monograph \cite{Chung_60}. Introduction of
taboo probabilities and hitting times under taboo provided a
powerful tool for study of functionals of Markov chains (see, e.g.,
\cite{Chung_60}, Ch.~2, Sec.~14), potential theory of Markov chains
(see, e.g., \cite{Kemeny_Snell_Knapp_Griffeath_76}, Ch.~4, Sec.~6),
trajectory properties (see, e.g., \cite{Zubkov_80}), matrix analytic
methods in stochastic modeling (see, e.g.,
\cite{Latouche_Ramaswami_99}, Ch.~3, Sec.~5), etc. As far as we know
the formula for probability of finiteness of a hitting time was
derived only for the cases of empty taboo set (see \cite{Chung_60},
Ch.~2, Sec.~12) and the taboo set consisting of a single state (see
\cite{B_SAB_12}). Now we complete the general picture. The results
are formulated as three theorems. The first one gives a
representation for the probability of finiteness of a hitting time
under taboo via taboo probabilities. The second theorem demonstrates
the relations between such probabilities with different initial and
target states when the taboo set changes by a single state. The
latter result allows to construct a finite iterative scheme to
evaluate the probability under consideration when either taboo set
or its complement are finite. Theorem \ref{T:transient} covers an
important particular case for a singleton taboo set. The proofs
involve the Laplace-Stieltjes transform of functions appearing in a
system of Chung's integral equations of convolution type.

Our interest in hitting times under taboo is motivated by their
application to effective classification (see \cite{Bul_Rimini_13})
of branching processes with finitely many catalysts. For \emph{a
single catalyst}, the model was described in \cite{DR_13}, although
in a more restrictive framework called \emph{branching random walk}
on $\mathbb{Z}^d$, $d\in\mathbb{N}$, it was proposed in
\cite{Yarovaya_91}. It turns out that in the branching random walk
with a single catalyst the number of particles at the point of
catalysis, coinciding with the process start point, can be
investigated by means of a due Bellman-Harris process with two types
of particles (see \cite{TV_SAM_13,VTY}). However, the study of other
local characteristics of the process with arbitrary start point can
be performed in terms of a Bellman-Harris process with six types of
particles (see \cite{B_Doklady_12,B_JTP_13}). For the employment of
such auxiliary processes, analysis of hitting times under taboo
\emph{for a random walk} on $\mathbb{Z}^d$ is indispensable as was
shown in \cite{B_SAB_12} and \cite{TV_SAM_13}. Note in passing that
the results announced in \cite{Bul_Rimini_13} enable us to
distinguish three types of asymptotic (as time tends to infinity)
behavior of particles population in \emph{branching processes with
finite set of catalysts} producing taboo sets. The type is
classified by the value of the Perron root (being less, equal or
greater than $1$) of a specified matrix with entries explicitly
depending on the studied probabilities of finiteness of hitting
times under taboo. So, intending to establish properties of
branching processes with several catalysts (and the particles
movement governed by a Markov chain), it is natural to prepare tools
by treatment of hitting times under taboo for a Markov chain.

\section{Main Results and Proofs}

We consider an irreducible continuous-time Markov chain
${\xi=\{\xi(t),t\geq0\}}$ generated by a conservative $Q$-matrix
$A=(a(x,y))_{x,y\in S}$ having finite negative diagonal elements.
Here $S$ is a finite or denumerable set. For $x\in S$, let
$\tau_x:=\mathbb{I}\{\xi(0)=x\}\inf\{t\geq0:\xi(t)\neq x\}$ where
$\mathbb{I}\{B\}$ stands for the indicator of a set $B$. The
stopping time $\tau_x$ (with respect to the natural filtration of
the process $\xi$) is \emph{the first exit time from $x$} and
$\mathbb{P}_x(\tau_x\leq t)=1-e^{a(x,x)t},$ $t\geq0$, (see, e.g.,
Theorem 5 in \cite{Chung_60}, Ch.~2, Sec.~5) where
$\mathbb{P}_x(\cdot):=\mathbb{P}(\cdot|\xi(0)=x)$. Following Chung's
notation in \cite{Chung_60}, Ch.~2, Sec.~11, for an arbitrary,
possibly empty set $H\subset S$ (``$\subset$'' always means
``$\subseteq$'') called henceforth the \emph{taboo set} and for
$t\geq0$, denote by
$$_H p_{x y}(t):=\mathbb{P}_x(\xi(t)=y,\;\xi(u)\notin H,\;min\{\tau_x,t\}<u<t),\quad x,y\in S,$$
the \emph{transition probability from $x$ to $y$ in time $t$ under
taboo $H$}. In the case $H=\varnothing$ the function $p_{x
y}(\cdot)={_\varnothing p}_{x y}(\cdot)$ is an ordinary transition
probability. Note that $_H p_{x y}(\cdot)\equiv0$ for $x\in S$,
$y\in H$ and $x\neq y$ whereas $_H p_{x x}(t)=e^{a(x,x)t}$ for $x\in
H$ and $_H p_{x x}(t)\geq e^{a(x,x)t}$ for $x\notin H$, $t\geq0$.
Set
$$_H\tau_{x y}:=\mathbb{I}\{\xi(0)=x\}\inf\{t\geq\tau_x:\xi(t)=y,\;\xi(u)\notin
H,\;\tau_x<u<t\},\; x,y\in S,$$ where, as usual, we assume that
$\inf\{t\in\varnothing\}=\infty$. The stopping time $_H\tau_{x y}$
is \emph{the first entrance time from $x$ to $y$ under taboo $H$}
whenever $x\neq y$ and is \emph{the first return time to $x$ under
taboo $H$} when $x=y$. Let ${_H F_{x y}(t):=\mathbb{P}_x(_H\tau_{x
y}\leq t)}$ and $F_{x y}(t):=\mathbb{P}_x({_\varnothing\tau}_{x
y}\leq t)$, $t\geq0$, be (improper) distribution functions of
$_H\tau_{x y}$ and ${_\varnothing\tau}_{x y}$, respectively.
Clearly, $_H\tau_{x y}={_{H\setminus\{y\}}\tau_{x y}}$ almost surely
(a.s.) for $y\in H$ and, consequently, ${_H F_{x
y}(\cdot)\equiv{_{H\setminus\{y\}}}F_{x y}(\cdot)}$ for $y\in H$.
So, it is sufficient to consider $_H F_{x y}(\cdot)$ for $x,y\in S$
and $y\notin H$.

According to Theorem 8 in \cite{Chung_60}, Ch.~2, Sec.~11, the
following \emph{first entrance formulae} are true for $x,y\in S$,
$z\notin H$ and $t\geq0$
\begin{eqnarray}
_H p_{x y}(t)&=&_{z,H}p_{x y}(t)+\int\nolimits_{0}^{t}{_H
p_{z y}(t-u)\,d\,_H F_{x z}(u)},\label{int_Hpxy(t)}\\
_H F_{x y}(t)&=&_{z,H}F_{x y}(t)+\int\nolimits_{0}^{t}{_H F_{z
y}(t-u)\,d\,_{y,H} F_{x z}(u)},\quad z\neq y,\label{int_HFxy(t)}
\end{eqnarray}
where we write $_{z,H}p_{x y}(t)$ instead of $_{z\cup H}p_{x y}(t)$
and similarly for other symbols. Prior to applying the Laplace
transform to functions from (\ref{int_Hpxy(t)}) and
(\ref{int_HFxy(t)}),  set
$$\widehat{p}(\lambda):=\int\limits_{0}^{\infty}{e^{-\lambda
t}p(t)\,dt},\;\;
\widehat{F}(\lambda):=\int\limits_{0-}^{\infty}{e^{-\lambda t}\,d
F(t)},\;\; P(t):=\int\limits_{0}^{t}{p(u)\,du},\;\;\lambda,t>0,$$
for any taboo probability $p$ and distribution function $F$.

Recall (see, e.g., \cite{Chung_60}, Ch.~2, Sec.~10) that the
irreducible Markov chain $\xi$ is recurrent (i.e.
$\mathbb{P}_x(\mbox{the set }\{t\geq\tau_x: \xi(t)=y\}\mbox{ is
unbounded})=1$ for any ${x,y\in S}$) iff $P_{x y}(\infty)=\infty$
where $P_{x y}(\infty):=\lim\nolimits_{t\to\infty}{P_{x y}(t)}$. In
a similar way, $\xi$ is transient (i.e. for each $x,y\in S$ one has
${\mathbb{P}_x(\mbox{the set }\{t\geq\tau_x: \xi(t)=y\}\mbox{ is
unbounded})=0}$) iff $P_{x y}(\infty)<\infty$. We stress that $\xi$
is either recurrent or transient (see, e.g., Theorem 4 and Corollary
2 in \cite{Chung_60}, Ch.~2, Sec.~10). For the properties of the
so-called \emph{Green function} ${G(x,y)\!:=P_{x y}(\infty)}$,
${x,y\in S}$, see, e.g., \cite{Lawler_Limic_10}, Ch.~4, Sec.~1
and~2. In accordance with Theorems~1,~3 and relation (5) in
\cite{Chung_60}, Ch.~2, Sec.~12, identity (\ref{int_Hpxy(t)})
implies that
\begin{eqnarray}
F_{x y}(\infty)=1& &\mbox{if}\quad\xi\;\mbox{is recurrent},\label{Fxy(infty)=recurrent}\\
F_{x y}(\infty)=\frac{G(x,y)}{G(y,y)}\in(0,1)&
&\mbox{if}\quad\xi\;\mbox{is
transient and}\; x\neq y,\quad\label{Fxy(infty)=transient}\\
F_{x x}(\infty)=1+\frac{1}{a(x,x)G(x,x)}\in(0,1)&
&\mbox{if}\quad\xi\;\mbox{is transient}.\label{Fxx(infty)=transient}
\end{eqnarray}

Our aim is to find the value $_H F_{x y}(\infty)$ being the
probability of finiteness of $_H\tau_{x y}$ conditioned on
$\{\xi(0)=x\}$. As was already noted it is sufficient to consider
$x,y\in S$ and $y\notin H$. The following statement shows that $_H
F_{x y}(\infty)$ can be expressed in terms of $_H P_{x y}(\infty)$
and $_H P_{y y}(\infty)$ similarly to (\ref{Fxy(infty)=transient})
and (\ref{Fxx(infty)=transient}).
\begin{thm}\label{T:H_Fxy(infty)=}
For any nonempty taboo set $H$ and $x,y\in S$, $y\notin H$, one has
\begin{eqnarray}
_H F_{x y}(\infty)=\frac{_H P_{x y}(\infty)}{_H
P_{y y}(\infty)},& &x\neq y,\label{T:_H_Fxy(infty)=}\\
_H F_{x x}(\infty)=1+\frac{1}{a(x,x){_H P_{x x}}(\infty)}\in[0,1),&
&x\notin H,\label{T:_H_Fxx(infty)=}
\end{eqnarray}
where $0\leq{_H P_{x y}}(\infty)<\infty$ and $0<{_H P_{y
y}}(\infty)<\infty$.
\end{thm}
\noindent{\sc Proof.} By Theorem 5 in \cite{Chung_60}, Ch.~2,
Sec.~11, the inequality ${_H P_{x y}(\infty)<\infty}$ is valid for
any \emph{nonempty} set $H$ and each $x,y\in S$. Setting $z=y$ we
apply the Laplace transform to both parts of (\ref{int_Hpxy(t)}).
Using the convolution property of the Laplace transform (see, e.g.,
\cite{Feller_71}, Ch.~13, Sec.~2, property (i)) we get
\begin{equation}\label{_H_widehat_F_x,y(lambda)=}
_H\widehat{F}_{x y}(\lambda)=\frac{_H\widehat{p}_{x
y}(\lambda)-{_{y,H}\widehat{p}_{x y}}(\lambda)} {_H\widehat{p}_{y
y}(\lambda)},\quad\lambda>0.
\end{equation}
Since $_{y,H}p_{x y}(\cdot)\equiv0$ for $x\neq y$, relation
(\ref{_H_widehat_F_x,y(lambda)=}) implies (\ref{T:_H_Fxy(infty)=})
due to identity
${F(\infty)=\lim\nolimits_{\lambda\to0+}{\widehat{F}(\lambda)}}$ for
a distribution function $F$ having support in $[0,\infty)$. The
identity holds by the monotone convergence theorem applied to the
Lebesgue integral representing $\widehat{F}(\lambda)$. Equality
(\ref{_H_widehat_F_x,y(lambda)=}) also entails
(\ref{T:_H_Fxx(infty)=}), since $_{x,H}p_{x x}(t)=e^{a(x,x)t}$,
$t\geq0$, and $_{x,H}\widehat{p}_{x
x}(\lambda)=(\lambda-a(x,x))^{-1}$, $\lambda\geq0$.
Theorem~\ref{T:H_Fxy(infty)=} is proved. $\square$

In the next theorem the value $_H F_{x y}(\infty)$ is expressed in
terms of $_{H'} F_{x' y'}(\infty)$ with appropriate choice of a
collection of states $x',y'\in S$ and a certain set $H'$ such that
$H'\subset H$ or $H\subset H'$. Thus, for a finite nonempty set $H$,
the evaluation of $_H F_{x y}(\infty)$ can be reduced to the case
when $H$ consists of a single state. The same procedure can be
performed when $S\setminus H$ is a finite set. Below we use the
Kronecker delta $\delta_{xy}$ for $x,y\in S$.
\begin{thm}\label{T:recursive}
If $H$ is a nonempty subset of $S$ and $x,y,z\in S$, $y,z\notin H$,
$z\neq y$, then
\begin{equation}\label{T:z,H_Fxy=}
_{z,H}F_{x y}(\infty)=\frac{_H F_{x y}(\infty)-{_H F_{x
z}(\infty)}{_H F_{z y}(\infty)}}{1-{_H F_{y z}(\infty)}{_H F_{z
y}(\infty)}}
\end{equation}
where ${_H F_{y z}(\infty)}{_H F_{z y}(\infty)}<1$. If $H$ is any
subset of $S$ and ${x,y\in S}$, $x\notin H$, $x\neq y$, then
\begin{equation}\label{T:H_Fxy=}
_H F_{x y}(\infty)=\frac{_{x,H}F_{x y}(\infty)} {1-{_{y,H}F_{x
x}}(\infty)}.
\end{equation}
Moreover, for any $H\subset S$ and $x,y\in S$ one has
\begin{equation}\label{T:H_Fxy=sum}
_H F_{x
y}(\infty)=(\delta_{xy}-1)\frac{a(x,y)}{a(x,x)}-\sum\limits_{z\in
S,\,z\neq x,\,z\neq y,\,z\notin H}{\frac{a(x,z)}{a(x,x)}{_H F_{z
y}}(\infty)}.
\end{equation}
\end{thm}
\noindent{\sc Proof.} Due to (\ref{int_HFxy(t)}) we have the
following system of two linear integral equations in functions
$_{z,H}F_{x y}(\cdot)$ and $_{y,H}F_{x z}(\cdot)$
$$
\left\{
\begin{array}{lcl}
_H F_{x y}(t)&=&_{z,H}F_{x y}(t)+\int\nolimits_{0}^{t}{_H F_{z y}(t-u)\,d\,_{y,H}F_{x z}(u)},\\
_H F_{x z}(t)&=&_{y,H} F_{x z}(t)+\int\nolimits_{0}^{t}{_H F_{y
z}(t-u)\,d\,_{z,H} F_{x y}(u)}.
\end{array}
\right.
$$
Applying the Laplace-Stieltjes transform and using its convolution
property (see, e.g., \cite{Feller_71}, Ch.~13, Sec.~2, property (i))
we get a new system of equations in $_{z,H}\widehat{F}_{x
y}(\lambda)$ and $_{y,H}\widehat{F}_{x z}(\lambda)$
$$
\left\{
\begin{array}{lcl}
_H\widehat{F}_{x y}(\lambda)&=&_{z,H}\widehat{F}_{x y}(\lambda)+
{_{y,H}\widehat{F}}_{x z}(\lambda)_H\widehat{F}_{z y}(\lambda),\\
_H\widehat{F}_{x z}(\lambda)&=&_{y,H}\widehat{F}_{x
z}(\lambda)+{_{z,H}\widehat{F}}_{x y}(\lambda)_H\widehat{F}_{y
z}(\lambda).
\end{array}
\right.
$$
Solving this system we obtain
$$_{z,H}\widehat{F}_{x y}(\lambda)=\frac{_H\widehat{F}_{x y}(\lambda)
-{_H\widehat{F}}_{x z}(\lambda)_H\widehat{F}_{z y}(\lambda)}{1
-{_H\widehat{F}}_{y z}(\lambda)_H\widehat{F}_{z y}(\lambda)}.$$
Letting $\lambda\to0+$ in the latter relation we come to
(\ref{T:z,H_Fxy=}). The inequality ${{_H F_{y z}(\infty)}{_H F_{z
y}(\infty)}<1}$ holds true, since ${_H F_{y z}(\infty)}{_H F_{z
y}(\infty)}\leq{_H F_{y y}(\infty)}$ in view of (\ref{int_HFxy(t)})
and $_H F_{y y}(\infty)<1$ by virtue of (\ref{T:_H_Fxx(infty)=}). We
come to relation (\ref{T:H_Fxy=}) applying the Laplace-Stieltjes
transform to (\ref{int_HFxy(t)}) when $z=x$ and letting
${\lambda\to0+}$. The claim (\ref{T:H_Fxy=sum}) ensues from the
identity
\begin{eqnarray*}
_H
F_{x y}(t)&=&(\delta_{xy}-1)\left(1-e^{a(x,x)t}\right)\frac{a(x,y)}{a(x,x)}\\
&-&\sum\limits_{z\in S,\,z\neq x,\,z\neq y,\,z\notin
H}{\frac{a(x,z)}{a(x,x)}\int\nolimits_{0}^{t}{\left(1-e^{a(x,x)(t-u)}\right)d{_H
F_{z y}}(u)}}
\end{eqnarray*}
that is valid due to the strong Markov property of $\xi$ (see, e.g.,
Theorem~3 in \cite{Chung_60}, Ch.~2, Sec.~9) involving the stopping
time $\tau_x$. The proof is complete. $\square$

The following result can be viewed as the complement to relation
(\ref{T:z,H_Fxy=}) when $H$ is an empty set and $\xi$ is a transient
Markov chain. The last hypothesis  permits us to obtain the formula
involving the Green functions.
\begin{thm}\label{T:transient}
Let $\xi$ be a transient Markov chain and $x,y,z\in S$. Then $_z
F_{x y}(\infty)\in[0,1)$ and
\begin{eqnarray}
_z F_{x
y}(\infty)=\frac{G(x,y)G(z,z)-G(x,z)G(z,y)}{G(z,z)G(y,y)-G(y,z)G(z,y)},\:
x\neq y,\:
x\neq z,&y\neq z,&\quad\label{T:zFxy=transient}\\
_z F_{y y}(\infty)=1+\frac{G(z,z)}{a(y,y)\left(G(y,y)G(z,z)-
G(y,z)G(z,y)\right)},&y\neq z,&\quad\label{T:zFyy=transient}\\
_z F_{z
y}(\infty)=-\frac{G(z,y)}{a(z,z)\left(G(y,y)G(z,z)-G(y,z)G(z,y)\right)},&y\neq
z.&\quad\label{T:zFzy=transient}
\end{eqnarray}
\end{thm}
\noindent{\sc Proof.} Examining the proof of Theorem
\ref{T:recursive}, we see that, for transient $\xi$, formula
(\ref{T:z,H_Fxy=}) is true even for $H=\varnothing$ as ${F_{y
z}(\infty)F_{z y}(\infty)<1}$ on account of
(\ref{Fxy(infty)=transient}). Thus, substituting
(\ref{Fxy(infty)=transient}) in (\ref{T:z,H_Fxy=}) we derive
(\ref{T:zFxy=transient}). According to (\ref{int_HFxy(t)}) and
(\ref{Fxy(infty)=transient}) one has $_z F_{x y}(\infty)\leq F_{x
y}(\infty)<1$. In a similar way we obtain (\ref{T:zFyy=transient})
and (\ref{T:zFzy=transient}) by employing
(\ref{Fxy(infty)=transient}) and (\ref{Fxx(infty)=transient}).
Theorem \ref{T:transient} is established. $\square$

For recurrent $\xi$, formula (\ref{T:z,H_Fxy=}) fails for
$H=\varnothing$, since ${F_{y z}(\infty)F_{z y}(\infty)=1}$ in view
of (\ref{Fxy(infty)=recurrent}). So, in general, there is no
counterpart of the previous theorem differing from assertion of
Theorem \ref{T:H_Fxy(infty)=} for a singleton taboo. However for
symmetric, space-homogeneous random walk on $\mathbb{Z}$ or
$\mathbb{Z}^2$ having finite variance of jump sizes (such random
walk is transient on $\mathbb{Z}^d$ with $d\geq3$) it is possible to
provide representation for $_z F_{x y}(\infty)$ alternative to those
in Theorem~\ref{T:H_Fxy(infty)=}. This follows from Theorems~1 and~2
in \cite{B_SAB_12}.

To conclude we return to the general case of Markov chains. In
contrast to $_H \tau_{x y}$ define the hitting time of state $y$
under taboo $H$ \emph{after the first exit out of the starting state
$x$} as
$$_H \overline{\tau}_{x y}:=\mathbb{I}\{\xi(0)=x\}\inf\{t\geq0:\xi(t+\tau_x)=y,\;\xi(u)\notin
H,\;\tau_x<u<t+\tau_x\}.$$ Such random variables arise naturally in
study of catalytic branching processes. Evidently, $_H\tau_{x
y}\!=\!\tau_{x}+{_H\overline{\tau}_{x y}}$ and
$\mathbb{P}_x({_H\overline{\tau}_{x
y}}=0)\!=\!(\delta_{xy}\!-\!1)a(x,y)a(x,x)^{\!-1}\!$. Moreover, by
virtue of the strong Markov property of $\xi$ random variables
$\tau_{x}$ and $_H\overline{\tau}_{x y}$ are independent. Therefore,
taking into account the convolution formula and the expression for
$\mathbb{P}_x(\tau_x\leq t)$ we get
\begin{equation}\label{convolution}
_H F_{x
y}(t)=\int\nolimits_{0-}^{t}{\left(1-e^{a(x,x)(t-u)}\right)d{_H\overline{F}_{x
y}(u)}}
\end{equation}
where ${_H\overline{F}_{x
y}(t):=\mathbb{P}_x\left(_H\overline{\tau}_{x y}\leq t\right)}$,
$t\geq0$. Hence, ${_H\overline{F}_{x y}(\infty)={_H F_{x
y}}(\infty)}$ for any ${x,y\in S}$, $H\subset S$, and the assertions
of Theorems \ref{T:H_Fxy(infty)=}~--~\ref{T:transient} hold true if
$_H F_{x y}(\infty)$ is replaced by $_H\overline{F}_{x y}(\infty)$.
Note also that on account of (\ref{convolution}) the distribution
function $_H F_{x y}(\cdot)$ has a bounded density. Consequently,
the function $_H\overline{F}_{x y}(\cdot)$ has also a density (which
is not bounded in general) in view of the following equality
$$_H
\overline{F}_{x y}(\infty)-{_H}\overline{F}_{x
y}(t)=\sum\limits_{z\in S,\;z\notin H,\;z\neq x,\;z\neq
y}{\frac{a(x,z)}{-a(x,x)}(_H F_{z y}(\infty)-{_H}F_{z y}(t))}$$
derived similarly to Lemma~3 in \cite{B_SAB_12}. Thus the results
established for $_H \tau_{x y}$ are valid for $_H\overline{\tau}_{x
y}$ as well.

\section{Acknowledgements}

The author is grateful to Professor V.A.Vatutin for useful
discussions and to the Reviewers for valuable remarks leading to
improvement of the text exposition.

\end{document}